\documentclass[12pt,reqno]{amsart}

\usepackage{amssymb}

\newcommand{\ds}{\displaystyle}     

\DeclareMathAlphabet{\E}{U}{eus}{m}{n}     

\newcommand{\PP}{{\mathbb P}}

\newcommand{\N}{{\mathbb N}}

\newcommand{\kk}{{\Bbbk}}

\newcommand{\Q}{{\mathfrak Q}}

\newcommand{\la}{\langle}
\newcommand{\ra}{\rangle}

\newtheorem{thm}{Theorem}[section]
\newtheorem{lemma}[thm]{Lemma}

\theoremstyle{definition}
\newtheorem{defn}[thm]{Definition}

\newtheorem{rmk}[thm]{Remark}

\newtheorem{example}[thm]{Example}

\newtheorem{Pf}{Proof$\!\!$}         
         
\newenvironment{pf}{\begin{Pf}}{\qed\end{Pf}}

\theoremstyle{remark}
\newtheorem{ack}{Acknowledgment\hspace*{-0.4em}}

\DeclareMathSymbol{\twoheadrightarrow}  {\mathrel}{AMSa}{"10}

\newcounter{letter}
\renewcommand{\theletter}{\rom{(}\alph{letter}\rom{)}}

\newcounter{rnum}
\renewcommand{\thernum}{\rom{(}\roman{rnum}\rom{)}}

\begin{document}
\baselineskip24pt


\title[GSCAs that are Twists of GCAs]%
{Graded Skew Clifford Algebras that are Twists\\[2mm] of Graded Clifford 
Algebras}

\subjclass[2010]{16S38, 16S37, 16S36}%
\keywords{Regular algebra, Clifford algebra, quadric, twist%
\rule[-3mm]{0cm}{0cm}}%

\maketitle

\vspace*{0.1in}

\baselineskip15pt

\renewcommand{\thefootnote}{\fnsymbol{footnote}}
\centerline{\sc Manizheh Nafari\footnote{Some of this research was carried
out while the first author worked at the University of Texas at
Arlington.}}
\centerline{Department of Mathematics and Statistics,}
\centerline{University of Toledo, Toledo, OH 43606-3390}
\centerline{{\sf manizheh.nafari@utoledo.edu}}

\bigskip
\centerline{and}
\bigskip

\centerline{\sc Michaela Vancliff\footnote{The second author was
supported in part by NSF grant DMS-0900239.}}
\centerline{Department of Mathematics, P.O.~Box 19408}
\centerline{University of Texas at Arlington,
Arlington, TX 76019-0408}
\centerline{{\sf vancliff@uta.edu \qquad www.uta.edu/math/vancliff}}

\setcounter{page}{1}
\thispagestyle{empty}

\bigskip
\bigskip

\begin{abstract}
\baselineskip14pt
In 2010, a quantized analog of a graded Clifford algebra (GCA), called a 
graded skew Clifford algebra (GSCA), was proposed by Cassidy and Vancliff,
and many properties of GCAs were found to have counterparts for GSCAs.
In particular, a GCA is a finite module over a certain commutative 
subalgebra~$C$, while a GSCA is a finite module over a (typically 
non-commutative) analogous subalgebra $R$. We consider the case that a 
regular GSCA is a twist of a GCA by an automorphism, and we 
prove, in this case, $R$ is a skew polynomial ring and a twist of 
$C$ by an automorphism. 
\end{abstract}

\bigskip
\bigskip

\baselineskip17.5pt


\section*{Introduction}
By 2011, in \cite{CV,NVZ}, it had been proved that almost all quadratic 
regular algebras of global dimension two or three may be classified using
certain non-commutative algebras called graded skew Clifford algebras 
(GSCAs). The latter algebras were first defined by Cassidy and
Vancliff in 
\cite{CV}, and may be viewed as a quantized analog of a graded Clifford 
algebra (GCA).  Many properties of GCAs were found in \cite{CV} to have 
counterparts for GSCAs; in particular, a GCA is quadratic and regular if 
and only if its associated quadric system has no base points, whereas a
GSCA is quadratic and regular if and only if its associated
(non-commutative) quadric system is normalizing and has no base points
(see Theorem~\ref{regular}).  Moreover, a GCA is a finite module over a 
certain commutative subalgebra~$C$, while a GSCA is a finite module over 
a (typically non-commutative) analogous
subalgebra $R$. In this article, we consider the case that a regular 
GSCA~$A$ is a twist by an automorphism of a GCA~$B$. In this
setting, we prove, in Theorem~\ref{mainthm}, that 
$R$~is a skew polynomial ring and is a twist of~$C$ by an automorphism.
We also demonstrate in Example~\ref{eg} that this can fail if 
$A$~is not a twist of a GCA.

This article consists of two sections: in Section 1, notation and 
terminology are defined, while Section 2 is devoted to proving our main 
result, which is given in Theorem~\ref{mainthm}.

\begin{ack} \  
The authors thank S.~P.~Smith of the University of Washington for
the suggestion to study the subalgebra $R$.
\end{ack}


\medskip

\section{Definitions}\label{Sec1}

In this section, we introduce the algebras to be discussed in the paper,
and some known results concerning them, including the connection between 
their homological properties and certain associated geometric data.

Throughout the article, $\kk$~denotes an algebraically closed field such
that char$(\kk)\neq~2$, and $M(n,\, \kk)$ denotes the vector space of 
$n \times n$ matrices with entries in $\kk$.  If $M$ is a matrix, then 
$M_{ij}$ will denote the $ij$'th entry of $M$.
For a graded 
$\kk$-algebra~$E$, the span of the homogeneous elements in $E$ of degree 
$i$ will be denoted $E_i$, and if $F$ is any ring or vector space, then 
$F^\times$ will denote the nonzero elements in $F$. 

\bigskip

Let $N_1 , \ldots , N_n \in M(n,\, \kk)$ denote symmetric matrices. By
definition (c.f., \cite{L}), a graded Clifford algebra (GCA) is the 
$\kk$-algebra $B$ on degree-one generators $X_1, \ldots , X_n$ and on 
degree-two generators $Y_1, \ldots , Y_n$ with defining relations given 
by
\begin{enumerate}
\item[(a)]
$\ds X_i X_j + X_j X_i = \sum_{k=1}^n (N_k)_{ij} \, Y_k$ for all $i, j =
1, \ldots , n$, and\\[-2mm]
\item[(b)] $Y_k$ central for all $k = 1, \ldots , n$. 
\end{enumerate}
We write $C$ for the subalgebra of $B$ generated by $Y_1, \ldots , Y_n$,
and note that $C$ is a polynomial ring and that $B$ is a finite 
module over $C$. Results on GCAs can be found in \cite{SV,VVW}.

The notion of graded skew Clifford algebra is similarly defined, 
but uses a generalization of the notion of symmetric matrix as 
follows. Let $\mu 
\in M(n,\, \kk)$, where $\mu_{ij}\mu_{ji} = 1$ for all $i,\ j$ with 
$i \neq j$. In \cite{CV}, a matrix $M \in M(n,\, \kk)$ is defined to be
$\mu$-symmetric if $M_{ij} = \mu_{ij}M_{ji}$ for all $i,\ j$.

\begin{defn}\cite{CV}\label{gsca} \  
Suppose, additionally, that 
$\mu_{ii} = 1$ for all $i$ and let $M_1, \ldots , M_n \in
M(n,\, \kk)$ be $\mu$-symmetric matrices. A graded skew Clifford algebra
(GSCA) is the $\kk$-algebra $A$ on degree-one generators $x_1, \ldots , 
x_n$ and on degree-two generators $y_1, \ldots , y_n$ with defining 
relations given by
\begin{enumerate}
\item[(a)]
degree-two relations:
$\ds x_i x_j + \mu_{ij} x_j x_i = \sum_{k=1}^n (M_k)_{ij} \, y_k$ for 
all $i, j = 1, \ldots , n$, and\\[-2mm]
\item[(b)] 
degree-three and degree-four relations that guarantee
the existence of a normalizing sequence $\{y_1', \ldots ,
y_n'\}$ that spans $\sum_{k=1}^n \kk y_k$. 
\end{enumerate}
We refer to the subalgebra of $A$ generated by $y_1, \ldots , y_n$ 
as $R$, and note that $A$ is a finite module over $R$.
\end{defn}

Clearly, a GCA is a special case of a GSCA. It is proved in \cite{CV}
that it is possible for all the $y_k$ to belong to $(A_1)^2$, and this 
happens if and only if $M_1, \ldots , M_n$ are linearly independent.
However, even in this case, $A$ need not be a quadratic algebra. 

One may associate to a symmetric matrix a quadratic form and 
hence a quadric in $\PP^{n-1}$. Similarly, as was shown in \cite{CV}, one
may associate a (non-commutative) quadratic form and a (non-commutative) 
quadric to a $\mu$-symmetric matrix as follows. Let $S$ denote the 
$\kk$-algebra
on generators $z_1, \ldots , z_n$ with defining relations $z_j z_i =
\mu_{ij} z_i z_j$ for all $i,\ j$, where $\mu_{ii} = 1$ for all $i$.
If $z = (z_1, \ldots , z_n)^T$ and if $M \in M(n,\, \kk)$ is a 
$\mu$-symmetric matrix, then the map $M \to z^T M z$ from the vector space of 
$\mu$-symmetric $n \times n$ matrices to $S_2$ is an isomorphism of vector
spaces (\cite{CV}).  As in \cite{CV}, we call the elements of 
$S_2$ quadratic forms. 
If $q \in S_2$, then the intersection in $\PP(S_1^*) \times 
\PP(S_1^*)$ of the zero locus of the defining relations of $S$ with the 
zero locus of $q$ is called the quadric associated to $q$. 
If $q \in S_2$ is normal in $S$, then its quadric
parametrizes those point modules over $S$ that are annihilated by $q$;
thus, this notion of quadric generalizes the commutative definition.

\begin{defn} \hfill\\
\indent (a) \cite{CV} \  
The span of quadratic forms $q_1, \ldots , q_m\in S_2$ will be called the 
{\em quadric system} associated to $q_1, \ldots , q_m$.
If a quadric system is given by a normalizing sequence in $S$, then it is
called a {\em normalizing quadric system}.\\
\indent (b) \cite{CVc} \  
We define a {\em left base point} of a quadric system $\{ \Q \}
\subset S_2$ to be any left base-point module over $S/\la \Q \ra$;
that is, to be any 1-critical graded left module over $S/\la \Q \ra$
that is generated by its homogeneous degree-zero elements and which
has Hilbert series $H(t) = c/(1-t)$, for some $c \in \N$.
We say a quadric system is {\em left base-point free} if it has no left base 
points. Similarly, for right base point, etc.
\end{defn}

If $S$ is commutative, then the notions of ``base point'' and ``base-point 
free'' agree with their commutative counterparts, since the only base-point 
modules in this case are point modules.

By \cite[Corollary~11]{CVc}, a normalizing quadric system $\Q$ is left
base-point free if and only if $\dim_\kk(S/\la \Q \ra) < \infty$. Hence, such 
a quadric
system is left base-point free if and only if it is right base-point free.
In particular, the adjectives ``right'' and ``left'' may be dropped when
referring to a {\em normalizing} quadric system being base-point free.

This geometric data associated to a GCA or GSCA has fundamental influence
on homological data of the algebra as demonstrated in the next result;
the reader is referred to \cite{AS,Lev} for definitions of the terms.

\begin{thm}\label{regular}\hfill
\begin{enumerate}
\item[(a)] \cite{Aubry.Lemaire,L} \  
The GCA $B$ is quadratic, Auslander-regular of global dimension~$n$ and 
satisfies the Cohen-Macaulay property with Hilbert series $1/(1-t)^n$ if 
and only if its associated (commutative) quadric system is base-point free; in 
this case, $B$~is AS-regular and is a noetherian domain.
\item[(b)] \cite{CV,CVc} \  
The GSCA $A$ is quadratic, Auslander-regular of global dimension $n$ and 
satisfies the Cohen-Macaulay property with Hilbert series $1/(1-t)^n$ if 
and only if its associated quadric system is normalizing and base-point 
free; in this case, $A$~is AS-regular, is a noetherian domain and is
unique up to isomorphism.
\end{enumerate}
\end{thm}

In the next section, we will consider a GSCA that is a twist (in the
sense of Definition~\ref{twistdef}) of a GCA by an automorphism. 
\begin{defn}\label{twistdef} \cite[\S8]{ATV2} \  
Let $D = \bigoplus_{k\geq 0} D_k$ be a graded $\kk$-algebra and let 
$\phi$ be a graded degree-zero automorphism of $D$.  The twist $D'$
of $D$ by $\phi$ is a graded $\kk$-algebra that is the vector space 
$\bigoplus_{k \geq 0} D_k$ with a
new multiplication $*$ defined as follows: if $a' \in D'_i = D_i,\ b' \in 
D'_j = D_j$, 
then $a' * b' = (a \phi^i(b))'$, where the right-hand side is computed using 
the original multiplication in $D$ and $a$ is the image of $a'$ in $D$, etc.
The twist of a quadratic algebra is again a quadratic algebra.
\end{defn}
\noindent
For $\phi$ and $a$ as in Definition~\ref{twistdef}, we will write 
$a^\phi$ for $\phi(a)$.

We close this section with a simple lemma concerning GCAs that will
be useful in the next section.
\begin{lemma}\label{central}
Let $B$ be a GCA as above. If $a,\ b \in B_1$, then $ab + ba$ is central 
in $B$.
\end{lemma}
\begin{pf}
The result is a consequence of $X_i X_j + X_j X_i$ being central in
$B$ for all $i, \ j$.
\end{pf}

\bigskip


\section{The Main Theorem}\label{Sec2}

In this section, we compare the subalgebras $R$ and $C$ that are
defined in Section~\ref{Sec1}. We prove in Theorem~\ref{mainthm} that if 
the GSCA~$A$ is a twist by an automorphism of a regular GCA~$B$, then 
$R$~is a twist of~$C$ by an automorphism and is a skew polynomial ring
(that is, $R$ is a domain on $n$~generators with $n(n-1)/2$~defining
relations that guarantee that each generator is normal in~$R$).

Not surprisingly, the algebra $R$ is not, in general, a skew polynomial ring 
nor a twist of a polynomial ring, and we first demonstrate this via a simple 
example.

\begin{example}\cite[\S3.2]{N} \label{eg} \  
Let $n = 3$, $\mu \in M(3,\, \kk)$ be as above and let 
\[
M_1 = \begin{bmatrix} 
& & \\[-4.7mm]
2&0&0\\[1mm] 0&0&0\\[1mm] 0&0&0\end{bmatrix},
\qquad
M_2 = \begin{bmatrix} 
& & \\[-4.7mm]
0&0&0\\[1mm] 0&2&0\\[1mm] 0&0&0\end{bmatrix},
\qquad
M_3 = \begin{bmatrix} 
& & \\[-4.7mm]
0&1&0\\[1mm] \frac12 &0&0\\[1mm] 0&0&2\end{bmatrix},
\]
where $\mu_{13} = \mu_{23} = 1$ and $\mu_{12}=2$. By
Theorem~\ref{regular}(b), the GSCA $A$ associated to this data is
quadratic and regular, and so is the $\kk$-algebra on $x_1$, $x_2$, 
$x_3$ with defining relations
\[
x_1 x_2 + 2 x_2 x_1 = x_3^2,\qquad
x_1 x_3 + x_3 x_1 = 0,\qquad
x_2 x_3 + x_3 x_2 = 0,
\]
where $y_i = x_i^2$ for $i = 1,\ 2,\ 3$. By Lemma~\ref{twist} (below), 
since $\mu_{13} \neq \mu_{12}\mu_{23}$, the associated algebra~$S$ is not
a twist of a polynomial ring and so, by \cite[Proposition~4.5]{CV}, $A$ 
is not a twist of a GCA by an automorphism. Moreover, $y_3$ is central 
in~$R$, but no other element in $\sum_{k=1}^3 y_k$ is normal in~$R$
(this can be seen by using a computer-algebra program such as 
W.~Schelter's
Affine program and noting that any normal element in $R$ would be
normal in $R/\la y_3 \ra$ in order to simplify the computations involved).
It follows that $R$ is not a skew polynomial ring.  Moreover, there
is an insufficient number of relations in low degree amongst the $y_k$ 
for $R$ to be a twist of a polynomial ring.
\end{example}

\begin{lemma}\label{twist}
Let $\mu \in M(n,\, \kk)$ and $S$ be as in Section~\ref{Sec1}. The
algebra $S$ is a twist of the polynomial ring $K = \kk [Z_1, \ldots
, Z_n]$ by a graded automorphism $\sigma\in$ Aut$(K)$ of degree zero if 
and only if $\mu_{ik} = \mu_{ij}\mu_{jk}$ for all $i, \ j, \ k$; in this 
case, $\sigma|_{K_1}$ is semisimple, and, for all $i,\ j$, we have
$\mu_{ij} = \rho_i/\rho_j$, where $\rho_i \in \kk^\times$ and 
$\sigma(Z_i) = \rho_i Z_i$ for all $i$.
\end{lemma}
\begin{pf}
The first part of the result follows from \cite{ATV1}, since 
$\mu_{ik} = \mu_{ij}\mu_{jk}$ for all $i, \ j, \ k$ if and only if the 
point scheme of $S$ (or the zero locus of the defining relations of $S$) 
is isomorphic to $\PP^{n-1}$, and the latter holds if and only if $S$ is 
a twist of the polynomial ring on $n$ variables by an automorphism.

For the second part of the result, suppose $S$ is a twist of the 
polynomial ring $K = \kk [Z_1, \ldots , Z_n]$ by a graded automorphism 
$\sigma\in$ Aut$(K)$ of degree zero, where we identify $z_i$ and
$Z_i$ for all $i$. From the relations in $S$, we have 
\[
Z_j Z_i^\sigma = \mu_{ij} Z_i Z_j^\sigma \quad \text{in\ } K \tag{$*$}
\]
for all $i,\ j$. However, $K$ is a commutative unique factorization
domain and deg$(Z_i) =1$ for all $i$, so $Z_i$ is irreducible in $K$.
Moreover, if $i \neq j$, then $Z_i \nmid Z_j$, so, by \thetag{$*$},
$Z_i | Z_i^\sigma$ for all $i$. Since $Z_i^\sigma$ has degree one, 
$Z_i^\sigma \in \kk^\times Z_i$ for all $i$. Hence, $\sigma|_{K_1}$
is semisimple. Writing $Z_i^\sigma = \rho_i Z_i$, where $\rho_i \in
\kk^\times$ for all $i$, and substituting into \thetag{$*$} completes 
the proof.
\end{pf}

\begin{rmk}\label{ssimple}
Suppose that $B$ is a regular GCA (in the sense of Theorem~\ref{regular})
and that $A$ is a GSCA that is a twist of $B$ by a graded automorphism 
$\tau \in$
Aut$(B)$ of degree zero. As was shown in Section~\ref{Sec1}, there
is a skew polynomial ring $S$ associated to $A$. By
\cite[Proposition~4.5]{CV}, since $A$ is a twist of $B$ by $\tau$,
there exists a choice for $S$ so that $S$ is a twist of the
polynomial ring $K = \kk [Z_1, \ldots , Z_n]$ by $\tau^{-1}$ and
conversely. By Lemma~\ref{twist}, $\tau|_{K_1}$ is semisimple and,
for each $i$, we have $\tau(Z_i) = \lambda_i Z_i$ for some
$\lambda_i \in \kk^\times$ and $\mu_{ij} = \lambda_j/\lambda_i$
for all $i,\ j$. (In the notation of Lemma~\ref{twist}, $\lambda_i = 
\rho_i^{-1}$ for all $i$, since $\tau = \sigma^{-1}$.)
\end{rmk}

\begin{thm}\label{mainthm}
Suppose that $A$ is a regular GSCA in the sense of
Theorem~\ref{regular}(b) and that $R$ is the subalgebra of $A$
generated by the $y_k$ as in Definition~\ref{gsca}. If $A$ is a
twist of a GCA $B$ by a graded automorphism $\tau\in$ Aut$(B)$ of 
degree zero, then $R$~is a twist of the analogous subalgebra $C$ of $B$ 
generated by the $Y_k$ and is a skew polynomial ring (that is, $R$~is a
domain that has exactly $n(n-1)/2$~defining relations that guarantee that 
each $y_k$ is normal in~$R$).
\end{thm}
\begin{pf}
By Remark~\ref{ssimple}, we may assume that $S$ is a twist of the 
polynomial ring $K = \kk [Z_1, \ldots , Z_n]$ by $\tau^{-1}$, and that 
$\tau(Z_i) = \lambda_i Z_i$ for some $\lambda_i \in \kk^\times$ and that 
$\mu_{ij} = \lambda_j/\lambda_i$ for all $i,\ j$. If necessary, by
Lemma~\ref{central}, we
may also re-choose the $X_k\in B_1$ so that $\{ X_1 , \ldots , X_n\}$
is dual to the basis $\{ Z_1 , \ldots , Z_n\}$ for $K_1$ and so that
the degree-two relations of $B$ still have the form given in
Section~\ref{Sec1} (although the symmetric matrices $N_1, \ldots , N_n$
might change). With this choice of bases, it follows that $X_i^\tau
= \lambda_i X_i$ for all $i$, and that the twist of $X_i$ is $x_i$.  
Hence, 
\[
x_i x_j + \mu_{ij} x_j x_i  = x_i x_j + (\lambda_j/\lambda_i) x_j x_i
\in\kk^\times (x_i^\tau x_j + x_j^\tau x_i )
\tag{$**$}
\]
for all $i,\ j$. We will prove that each $x_i^\tau x_j + x_j^\tau
x_i = r_{ij}$  is a normal element of $A$. By Definition~\ref{gsca}
and \thetag{$**$},
$r_{ij} \in R$ for all $i,\ j$, so the subalgebra of $A$ generated by the
$r_{ij}$ is contained in $R$. Since $A$ is quadratic, each $y_k$ is
a function of the $r_{ij}$ and so $R$ is the subalgebra of $A$ generated 
by the $r_{ij}$. Moreover,  for all $i,\ j,\ k$, we have 
\[
\begin{array}{ll}
x_k r_{ij} & = x_k (x_i^\tau x_j + x_j^\tau x_i)\\[3mm]
& = X_k (X_i^{\tau^2} X_j^{\tau^2} + X_j^{\tau^2} X_i^{\tau^2})\\[3mm]
& = \lambda_i^2 \lambda_j^2 X_k (X_i X_j + X_j X_i)\\[3mm]
& = \lambda_i^2 \lambda_j^2 (X_i X_j + X_j X_i) X_k\\[3mm]
& = \lambda_k^{-2} \lambda_i^2 \lambda_j^2 
              (x_i^\tau x_j + x_j^\tau x_i) x_k\\[3mm]
& = \mu_{ki} \mu_{kj} r_{ij} x_k,
\end{array}
\]
so that the $r_{ij}$ are normal in $A$. It follows that
\[
r_{ij} r_{kp} = \mu_{ik} \mu_{jk} \mu_{ip} \mu_{jp} r_{kp} r_{ij}
              = \mu_{ik}^2 \mu_{jp}^2 r_{kp} r_{ij},
	      \tag{$\dag$}
\]
since $\mu_{ij} = \lambda_j/\lambda_i$ for all $i,\ j$. Thus, $R$ is
a skew polynomial ring. For all $i,\ j,\ k,\ p$, let $\nu_{ijkp} =
\mu_{ik}^2 \mu_{jp}^2$, so $\nu_{ijkp} \nu_{kpab}  = \nu_{ijab}$,
for all $i,\ j,\ k,\ p,\ a,\ b$. By Lemma~\ref{twist}, it follows that 
$R$ is a twist of a polynomial ring.

For all $i,\ j$, let $c_{ij}\in B$ denote the element that twists
to $r_{ij} \in A$; that is, 
\[
c_{ij} = X_i^\tau X_j^\tau + X_j^\tau X_i^\tau
= \tau(X_i X_j + X_j X_i)
\in \kk^\times(X_i X_j + X_j X_i)
\subset C.
\]
Moreover, since $B$ is quadratic, each $Y_k$ is a function of the 
$X_i X_j + X_j X_i$ and so a function of the $c_{ij}$. It follows
that $C$ is the subalgebra of $B$ generated by the $c_{ij}$, and so
$R$ is a twist of~$C$. By \thetag{$\dag$}, we have 
\[
c_{ij} c_{kp}^{\tau^2} = \nu_{ijkp} c_{kp} c_{ij}^{\tau^2},
\]
for all $i,\ j,\ k,\ p$. Defining $\tau' \in$ Aut$(C)$ by 
$\tau' (c_{ij}) = \lambda_i^2 \lambda_j^2 c_{ij}$, for all $i,\ j$ (so
$\tau' = \tau^2|_C$), we find that $R$ is a twist of $C$ by $\tau'$.
\end{pf}


\end{document}